\documentclass{article}
\begin{document}
\newtheorem{proposition}{Proposition}[section]
\newtheorem{definition}{Definition}[section]
\newtheorem{lemma}{Lemma}[section]
\newcommand{\xl}{\stackrel{\rightharpoonup}{\cdot}}
\newcommand{\xr}{\stackrel{\leftharpoonup}{\cdot}}
\newcommand{\xlplus}{\stackrel{\rightharpoonup}{+}}
\newcommand{\xrplus}{\stackrel{\leftharpoonup}{+}}
\newcommand{\Ll}{\stackrel{\rightharpoonup}{L}}
\newcommand{\Lr}{\stackrel{\leftharpoonup}{L}}
\newcommand{\Rl}{\stackrel{\rightharpoonup}{R}}
\newcommand{\Rr}{\stackrel{\leftharpoonup}{R}}

\title{\bf A class of group-like objects}
\author{Keqin Liu\\Department of Mathematics\\The
University of British Columbia\\Vancouver, BC\\
Canada, V6T 1Z2}
\date{October 16, 2003}
\maketitle

\begin{abstract} We introduce a class of group-like
objects and prove that Cayley Theorem on groups has a
counterpart in the class of group-like objects.
\end{abstract}

In the study of generalizing the Lie correspondence
between linear Lie groups and linear Lie algebras, we
have found a generalization of the notion of inverse
in group theory. The generalization of the notion of
inverse in group theory not only is indispensable in
our method of generalizing the Lie correspondence
between linear Lie groups and linear Lie algebras, but
also plays a key role in our study of group-like,
ring-like and field-like objects. In this paper, we
use the generalization of the notion of inverse to
introduce a class of group-like objects. The class of
group-like objects share many properties with groups.
The main result of this paper is that Cayley Theorem
on groups has a counterpart in the class of group-like
objects.

\section{Basic Definitions}

We first give the definition of the class of
group-like objects.

\begin{definition}\label{def0.1.1} Let $G$ be a
nonempty set together with two binary operations $\xl$
and $\xr$ on $G$. $G$ is called a {\bf digroup} with
the {\bf identity} $e$ under the two binary operations
if $e\in G$  and the following three properties are
satisfied

1. The two operations $\xl$ and $\xr$ are
diassociative, that is,
\begin{equation}\label{eq0.1.1} x\xl (y\xl z)=(x\xl
y)\xl z=x\xl (y\xr z)\end{equation}
\begin{equation}\label{eq0.1.2} (x\xr y)\xl z=x\xr
(y\xl z)\end{equation}
\begin{equation}\label{eq0.1.3} (x\xl y)\xr z)=(x\xr
y)\xr z=x\xr (y\xr z)\end{equation}
for all $x,y,z\in G$.

2. For all $x\in G$, we have
\begin{equation}\label{eq0.1.4} x\xl e =x=e\xr
x\end{equation}
\begin{equation}\label{eq0.1.5} x\xr e=e\xl
x\end{equation}

3. For each element $x$ in $G$, there is an element 
$x^{\check{-}1}$ in $G$, call the {\bf Liu inverse} of
$x$, such that
\begin{equation}\label{eq0.1.6} x^{\check{-}1}\xl x
=e= x\xr x^{\check{-}1}\end{equation}.
\end{definition}

The binary operations $\xl$ and $\xr$ are called the
{\bf left product} and the {\bf right product},
respectively. A digroup $G$ is also denoted by $(G,
e)$ or $(G, \xl ,  \xr , e)$.
A digroup $G$ is called a {\bf finite digroup} (or an
{\bf infinite digroup}) if $G$ is a finite (or
infinite) set. The cardinal number of a set will be
called the {\bf order} of the set.

\medskip
(\ref{eq0.1.1}), (\ref{eq0.1.2}) and (\ref{eq0.1.3})
consist of the  {\bf diassociative law}, which was
introduced by J.-L.Loday to study Leibniz
algebras(\cite{Loday1}). A set with two binary
operations $\xl $ and $\xr $ is called a {\bf
dimonoid} if the two binary operations 
$\xl $ and $\xr $ satisfy the  diassociative law
(\cite{Loday2}). An element $e$ of a dimonoid
satisfying (\ref{eq0.1.4}) is called a {\bf bar-unit}
(\cite{Loday2}).

\medskip
By Definition~\ref{def0.1.1}, a digroup $G$ is defined
with respect to an fixed element $e$ called the
identity of $G$. An Liu inverse of an element of $G$
is defined with respect to the identity $e$. The
following proposition establishes the uniqueness of an
Liu inverse.
\begin{proposition}\label{pr0.1.1} If $x$ is an
element of a digroup $(G., e)$, then the Liu inverse
of $x$ is unique. In other words, if $y_i\xl x=e=x\xr
y_i$ for $y_i\in G$ and $i=1, 2$, then $y_1=y_2$.
\end{proposition}

\medskip
\noindent
{\bf Proof} This is a direct consequence of
Definition~\ref{def0.1.1}. \qquad $\|$

\medskip
It is clear that a digroup $(G, \xl , \xr , e)$
becomes a group if $\xl =\xr$. The inverse of an
element $\alpha $ of a group will be called the {\bf
group inverse} of $\alpha$, and the ordinary notation
$\alpha ^{-1}$ will be used to denote the group
inverse of $\alpha $.

\begin{definition}\label{def0.1.2} Let $(G, \xl , \xr
, e)$ be a digroup.
\begin{description}
\item[(i)] Two elements $x$ and $y$ of $G$ are said to
be {\bf commutative} if
$$ x\xl y=y\xr x.$$ 
\item[(ii)] $G$ is said to be a {\bf commutative
digroup} if any two elements of $G$ are commutative.
\end{description}
\end{definition}

\medskip
If a digroup $G$ is commutative, then the additive
notations $\xlplus$ and $\xrplus$ are also used to
denote the two binary operations on $G$. Also in this
additive notations, we write $0$ for the identity, and
we write $\check{-}x$ for the Liu inverse of $x$ in
$G$.

\medskip
\noindent
{\bf Example 1} Let $M: =\{\, 0, \, a\,\}$ be a set of
two distinct elements. We define two binary operations
$\xlplus$ and $\xrplus$ on $M$ as follows:
$$
\begin{tabular}{|c||c|c|}\hline
$\xlplus$ &0&a\\
\hline\hline
0&0&0\\
\hline
a&a&a\\
\hline
\end{tabular}
\qquad\qquad\qquad 
\begin{tabular}{|c||c|c|}\hline
$\xrplus$ &0&a\\
\hline\hline
0&0&a\\
\hline
a&0&a\\
\hline
\end{tabular}
$$

It is easy to check that $M$ is a commutative digroup
with the identity $0$, and $0$ is the Liu inverse for
both $0$ and $a$. Since a digroup consisting of a
single element must be the identity group, $M$ is the
smallest digroup which is not a group. $\qquad \|$

\medskip
\noindent
{\bf Example 2} Let $N: =\{\, e, \, \alpha ,\, \beta
,\, \gamma ,\, \delta ,\, \epsilon \,\}$ be a set of
six distinct elements. We define two binary operations
$\xl$ and $\xr$ on $N$ as follows:
$$
\begin{tabular}{|c||c|c|c|c|c|c|}\hline
$\xl$
&e&$\alpha$&$\beta$&$\gamma$&$\delta$&$\epsilon$\\
\hline\hline
e&e&$\alpha$&$\alpha$&$\alpha$&e&e\\
\hline
$\alpha$&$\alpha$&e&e&e&$\alpha$&$\alpha$\\
\hline
$\beta$&$\beta$&$\delta$&$\delta$&$\delta$&$\beta$&$\beta$\\
\hline
$\gamma$&$\gamma$&$\epsilon$&$\epsilon$&$\epsilon$&$\gamma$&$\gamma$\\
\hline
$\delta$&$\delta$&$\beta$&$\beta$&$\beta$&$\delta$&$\delta$\\
\hline
$\epsilon$&$\epsilon$&$\gamma$&$\gamma$&$\gamma$&$\epsilon$&$\epsilon$\\
\hline
\end{tabular}
\qquad\qquad\qquad
\begin{tabular}{|c||c|c|c|c|c|c|}\hline
$\xr$
&e&$\alpha$&$\beta$&$\gamma$&$\delta$&$\epsilon$\\
\hline\hline
e&e&$\alpha$&$\beta$&$\gamma$&$\delta$&$\epsilon$\\
\hline
$\alpha$&$\alpha$&e&$\epsilon$&$\delta$&$\gamma$&$\beta$\\
\hline
$\beta$&$\alpha$&e&$\epsilon$&$\delta$&$\gamma$&$\beta$\\
\hline
$\gamma$&$\alpha$&e&$\epsilon$&$\delta$&$\gamma$&$\beta$\\
\hline
$\delta$&e&$\alpha$&$\beta$&$\gamma$&$\delta$&$\epsilon$\\
\hline
$\epsilon$&e&$\alpha$&$\beta$&$\gamma$&$\delta$&$\epsilon$\\
\hline
\end{tabular}
$$

One can check that $N$ is a digroup with the identity
$e$. Since 
$$\beta \xl \beta=\delta\ne \epsilon=\beta \xr \beta
,$$
$N$ is not commutative. By the way, we indicate that
$N$ is the unique non-commutative digroup with order
$6$, and the number {\bf 6}, which is the smallest
perfect number, is the smallest order among the orders
of non-commutative digroups. Hence, $N$ is the
smallest non-commutative digroup which is not a group.
$\qquad \|$

\begin{definition}\label{def1.2.1} Let $H$ be a
nonempty subset of a digroup $(G, \xl , \xr , e)$. If
$e\in H$ and $(H, \xl , \xr , e)$ is itself a digroup
with the identity e, we say that $H$ is a {\bf
subdigroup} of $G$. $H\le G$ or $G\ge H$ is used to
signify that $H$ is a subdigroup of $G$.
\end{definition}

\begin{proposition}\label{pr1.2.1} Let $H$ be a subset
of a digroup $(G, \xl , \xr , e)$. The following are
equivalent.
\begin{description}
\item[(i)] $H$ is a subdigroup of $G$.
\item[(ii)] $e\in H$ and $(H\xl H^{\check{-}1})\cup
(H^{\check{-}1}\xr H)\subseteq H$. 
\item[(iii)] $H$ is nonempty and $(H\ast H)\cup
H^{\check{-}1}\subseteq H$, where $\ast =\xl$ or
$\xr$.
\end{description}
\end{proposition}

\medskip
\noindent
{\bf Proof} (i)$\Rightarrow$ (ii): This is  clear by
Definition~\ref{def1.2.1}.

\medskip
(ii)$\Rightarrow$ (iii): For $h\in H$, we have 
$h^{\check{-}1}=
e\xl h^{\check{-}1}\in H\xl H^{\check{-}1}\subseteq
H$. Hence, we get 
$H^{\check{-}1}\subseteq H$.

For $h_1$, $h_2\in H$, we have
\begin{eqnarray*}
&&h_1\xl h_2=h_1\xl (e\xr h_2)=h_1\xl (e\xl h_2)=
h_1\xl (h_2^{\check{-}1})^{\check{-}1}\in H,
\end{eqnarray*}
which implies that $H\xl H\subseteq H$. Similarly, we
have $H\xr H\subseteq H$. 

This proves that (iii) is true.

\medskip
(iii)$\Rightarrow$ (i): Since $H\ne\emptyset$, there
is an element $h\in H$. Hence, $h^{\check{-}1}\in H$.
It follows that $e=h^{\check{-}1}\xl h\in H\xl
H\subseteq H$. By $H\ast H\subseteq H$, $H$ is closed
under both the left product $\xl$ and the right
product $\xr$. This proves that $(H, \xl , \xr , e)$
is a digroup. \qquad $\|$

\section{The counterpart of Cayley Theorem}

Let $\Omega$ be a set. The set of all maps from
$\Omega$ to $\Omega$ is denoted by
$\mathcal{T}(\Omega)$. Thus
$$ \mathcal{T}(\Omega):=\{\; f \; | \; \mbox{$f:
\Omega\to\Omega$ is a map} \; \}.$$

It is well known that $\mathcal{T}(\Omega)$ is a
semigroup with the identity $1$ under the product
$fg:=f\cdot g$, where $1$ is the identity map, and the
product $f\cdot g$ is defined by
$$(f\cdot g)(x): =f(g(x)) \quad\mbox{for $x\in
\Omega$}.$$ 

\medskip
Let $(G, \xl , \xr , e)$ be a digroup. For $a\in G$,
we define two maps $\Ll_a$ and $\Lr_a$ as follows:
$$ \Ll_a(x):=a\xl x, \quad \Lr_a(x):=a\xr x
\quad\mbox{for all $x\in G$}.$$
$\Ll_a$ and $\Lr_a$ are called the {\bf left
translations} determined by $a$. Let 
$$\Lr_G:=\{\; \Lr_a \; | \; a\in G \; \}, \quad
\Ll_G:=\{\; \Ll_a \; | \; a\in G \; \}.$$
Then both $\Lr_G$ and $\Ll_G$ are subsets of
$\mathcal{T}(G)$. 

\begin{proposition}\label{pr1.4.1} Let $(G, \xl , \xr
, e)$ be a digroup. If $a,b\in G$, then 
\begin{equation}\label{eq1.4.1} 
\Lr_{a\ast b}=\Lr_a\Lr_b, \quad \Ll_{a\xl
b}=\Ll_a\Ll_b=\Ll_a\Lr_b, \quad
\Lr_{a\xr b}=\Lr_a\Ll_b, \quad  
\end{equation}
\begin{equation}\label{eq1.4.2} 
\Lr_e=1, \quad \Ll_a\Ll_e=\Ll_a, \quad
\Ll_e\Ll_a=\Lr_a\Ll_e, \quad  
\end{equation}
\begin{equation}\label{eq1.4.3} 
\Lr_{a^{\check{-}1}}\Lr_a=1=\Lr_a\Lr_{a^{\check{-}1}},
\quad 
\Ll_{a^{\check{-}1}}\Ll_a=\Ll_e=\Lr_a\Ll_{a^{\check{-}1}}.
\end{equation}
where $\ast=\mbox{$\xl$ or $\xr$.}$
\end{proposition}

\medskip
\noindent
{\bf Proof} Each of  the equations above follows from
the definition of a digroup. For example, let us prove
that $\Ll_e\Ll_a=\Lr_a\Ll_e$. For $x\in G$, we have
\begin{eqnarray*}
&&(\Ll_e\Ll_a)(x)=\Ll_e(\Ll_a(x))=e\xl (a\xl x)
=(a\xl x)\xr e\\
&=&a\xr (x\xr e)
=a\xr (e\xl x)
=a\xr (\Ll_e(x))=\Lr_a(\Ll_e(x))=(\Lr_a\Ll_e)(x), 
\end{eqnarray*}
which implies that $\Ll_e\Ll_a=\Lr_a\Ll_e$. $\qquad
\|$

\medskip
The following three propositions summerize the main
properties of $\Ll_G$ and $\Lr_G$.

\begin{proposition}\label{pr1.4.2} If $(G, \xl , \xr ,
e)$ is a digroup, then $\Lr_G$ is a subgroup of
$\mathcal{T}(G)$ such that the group inverse of
$\Lr_a$ in the group $\Lr_G$ is
$\Lr_{a^{\check{-}1}}$, where $a\in G$ and
$a^{\check{-}1}$ is the Liu inverse of $a$ in the
digroup $G$.
\end{proposition}

\medskip
\noindent
{\bf Proof} By (\ref{eq1.4.1}), $\Lr_a\Lr_b=\Lr_{a\ast
b}\in\Lr_G$ for $a,b\in G$. Hence, $\Lr_G$ is closed
under the product in $\mathcal{T}(G)$. By
(\ref{eq1.4.3}), $\Lr_{a^{\check{-}1}}$ is the group
inverse of $\Lr_a$. Thus $\Lr_G$ is a subgroup of
$\mathcal{T}(G)$. $\qquad \|$

\begin{proposition}\label{pr1.4.3} If $(G, \xl , \xr ,
e)$ is a digroup, then $\Ll_G$ is a subsemigroup of
$\mathcal{T}(G)$ such that 
\begin{description}
\item[(i)] $\Ll_e$ is a right unit of $\Ll_G$, i.e.,
$$ \Ll_a\Ll_e=\Ll_a \quad\mbox{ for $a\in G$,}$$
\item[(ii)] $\Ll_{a^{\check{-}1}}$ is a left inverse
of $\Ll_a$ with respect to $\Ll_e$, i.e.,
$$\Ll_{a^{\check{-}1}}\Ll_a=\Ll_e.$$
\end{description}
\end{proposition}

\medskip
\noindent
{\bf Proof} For $a,b\in G$, we have
$\Ll_a\Ll_b=\Ll_{a\xl b}\in \Ll_G$ by (\ref{eq1.4.1}).
Hence, $\Ll_G$ is a subsemigroup of $\mathcal{T}(G)$.
$(i)$ and $(ii)$ follow from $(\ref{eq1.4.2})$ and
$(\ref{eq1.4.3})$, respectively.  $\qquad \|$

\begin{proposition}\label{pr1.4.4} If $(G, \xl , \xr ,
e)$ is a digroup, then the map $\phi: \Lr_G\to \Ll_G$
defined by
$$ \phi (\Ll_a):=\Lr_a \quad\mbox{for $a\in G$} $$
is a semigroup homomorphism satisfying
\begin{equation}\label{eq1.4.4} \phi (\Ll_e)=1, \quad
\Ll_e\Ll_a=\phi(\Ll_a)\Ll_e,\end{equation}
\begin{equation}\label{eq1.4.5}
\phi(\Ll_a)\Ll_{a^{\check{-}1}}=\Ll_e,\end{equation}
\begin{equation}\label{eq1.4.6} \Ll_a\phi
(\Ll_b)=\Ll_a\Ll_b, \quad
\phi(\phi(\Ll_a)\Ll_b)=\phi(\Ll_a)\phi(\Ll_b),\end{equation}
where $a, b\in G$.
\end{proposition}

\medskip
\noindent
{\bf Proof}   $\phi $ is well-defined. In fact, if
$\Ll_a=\Ll_b$ for $a, b\in G$, then $a=a\xl
e=\Ll_a(e)=\Ll_b(e)=b\xl e=b$. Hence,
$\phi(\Ll_a)=\phi(\Ll_b)$.

By (\ref{eq1.4.1}), $\phi(\Ll_a\Ll_b)=\phi(\Ll_{a\xl
b})=\Ll_{a\xr b}=\Lr_a\Lr_b=\phi(\Ll_a)\phi(\Ll_b)$.
Thus, $\phi$ is a semigroup homomorphism .

$\phi (\Ll_e)=\Lr_e=1$ by (\ref{eq1.4.2}). Also, we
have $\Ll_e\Ll_a=
\Ll_{e\xl a}=\Ll_{a\xr
e}=\Lr_a\Ll_e=\phi(\Ll_a)\Ll_e.$ This proves
(\ref{eq1.4.4}).

Since
$\phi(\Ll_a)\Ll_{a^{\check{-}1}}=\Lr_a\Ll_{a^{\check{-}1}}
=\Ll_e$, (\ref{eq1.4.5}) holds.

Finally, by (\ref{eq1.4.1}), we have 
$$ \Ll_a\phi (\Ll_b)=\Ll_a\Lr_b
=\Ll_{a\xl b}=\Ll_a\Ll_b$$
and
$$\phi(\phi(\Ll_a)\Ll_b)=\phi(\Lr_a\Ll_b)=\phi(\Ll_{a\xr
b})=\Lr_{a\xr b}=\Lr_a\Lr_b
=\phi(\Ll_a)\phi(\Ll_b),$$
which prove that (\ref{eq1.4.6}) holds.  $\qquad \|$

\bigskip
Let $(G, \xl , \xr , e)$ be a digroup. We define two
binary operations $\xl$ and $\xr$ on the set
$\Lr_G\times \Ll_G$ by
\begin{eqnarray}
\label{eq1.4.7} (\Lr_a, \Ll_b)\xl (\Lr_c,
\Ll_d)&=&(\Lr_{a\xl c}, \; \Ll_{b\xl d}),\\
\label{eq1.4.8} (\Lr_a, \Ll_b)\xr (\Lr_c,
\Ll_d)&=&(\Lr_{a\xr c}, \; \Ll_{b\xr d}),
\end{eqnarray}
where $a$, $b$, $c$, $d\in G$. 

It is clear that under the binary operations above,
$\Lr_G\times \Ll_G$ becomes a digroup with the
identity $(\Lr_e, \Ll_e)=(1, \Ll_e)$, and the Liu
inverse of $(\Lr_a, \Ll_b)$ is 
$(\Lr_{a^{\check{-}1}}, \Ll_{b^{\check{-}1}})$.

\medskip
The digroup $\Lr_G\times \Ll_G$ obtained from a
digroup $G$ has a natural action on the set $G\times
G$:
$$(\Lr_a, \Ll_b)(x,y):=(\Lr_a(x),
\Ll_b(y))\quad\mbox{for $(x,y)\in G\times G$}.$$
Hence, $\Lr_G\times \Ll_G$ can be regarded as a subset
of 
$$\mathcal{T}(G\times G):=\{\; f \; | \; \mbox{$f:
G\times G\to G\times G$ is a map} \; \}.$$
Note that the identity $(1, \Ll_e)$ of the digroup
$\Lr_G\times \Ll_G$ does not act like the identity map
on the set $G\times G$.

\begin{definition}\label{def1.4.1} A {\bf
homomorphism} $\eta$ from a digroup $(G, e_G)$ to a
digroup $(H, e_H)$ is a map from $G$ to $H$ such that 
\begin{eqnarray}
\eta(e_G)&=&e_H,\nonumber\\
\eta(x\ast y)&=&\eta(x)\ast\eta(y),\nonumber
\end{eqnarray}
where $x,y\in G$ and $\ast =\xl , \; \xr$. We say that
two digroups are {\bf isomorphic} if there is a
bijective homomorphism from one to the
other.\end{definition}

The following proposition gives a counterpart of
Cayley Theorem.

\begin{proposition}\label{pr1.4.5} Every digroup $(G,
e)$ is isomorphic a subdigroup of $\Lr_G\times \Ll_G$.
\end{proposition}

\medskip
\noindent
{\bf Proof} Let
$$\mathcal{D}(\Lr_G\times \Ll_G):=\{\; (\Lr_a, \Ll_a)
\; | \; a\in G \;\}.$$

By (\ref{eq1.4.7}) and (\ref{eq1.4.8}),
$\mathcal{D}(\Lr_G\times \Ll_G)$ is a subdigroup of
$\Lr_G\times \Ll_G$.

\medskip
Define $\eta : G\to \mathcal{D}(\Lr_G\times \Ll_G)$ by
$$\eta (a):=(\Lr_a, \Ll_a)\quad a\in G.$$

First, it is clear that $\eta(e)=(\Lr_e, \Ll_e)=(1,
\Ll_e)$ is the identity of the digroup
$\mathcal{D}(\Lr_G\times \Ll_G)$, and $\eta $ is
surjective.

Next, we have
\begin{eqnarray}
&&\eta(a\ast b)=(\Lr_{a\ast b}, \; \Ll_{a\ast
b})=(\Lr_a, \Ll_a)\ast (\Lr_b,
\Ll_b)=\eta(a)\ast\eta(b).\nonumber
\end{eqnarray}
Hence, $\eta$ preserves the binary operations.

Finally, $\eta(a)=\eta(b)\Rightarrow
\Ll_a=\Ll_b\Rightarrow a=b$. Thus, $\eta$ is
injective.

This proves that $\eta$ is a isomorphism from $G$ to
$\mathcal{D}(\Lr_G\times \Ll_G)$. 
$\qquad \|$

\medskip
If $a$ is an element of a digroup $(G, \xl , \xr ,
e)$, we may also define two maps $\Rl_a$ and
$\Rr_a\in\mathcal{T}(G)$ as follows:
$$ \Rl_a(x):=x\xl a, \quad \Rr_a(x):=x\xr a
\quad\mbox{for all $x\in G$}.$$
$\Rl_a$ and $\Rr_a$ are called the {\bf right
translations} determined by $a$. Let $\Rl_G$ and
$\Rr_G$ be two subsets of $\mathcal{T}(G)$ defined by
$$\Rl_G:=\{\; \Rl_a \; | \; a\in G \; \}, \quad
\Rr_G:=\{\; \Rr_a \; | \; a\in G \; \}.$$

A similar argument will show that $\Rr_G\times \Rl_G$
can be made into a digroup such that $G$ is isomorphic
to a subdigroup of $\Rr_G\times \Rl_G$.

\medskip
We finish this section with a construction of the
digrouop which generalizes the construction of the
digrouop  $\Lr_G\times \Ll_G$.

\medskip
Let $\Omega$ be a set. A triple $(\mathcal{S},
\mathcal{G}, \phi )$ is called a {\bf standard triple}
on $\Omega$ if it satisfies the following three
conditions:
\begin{enumerate}
\item $\mathcal{G}$ is a subgroup of
$\mathcal{T}(\Omega)$ with the identity $1$.
\item $\mathcal{S}$ is a subsemigroup of
$\mathcal{T}(\Omega)$ such that
\begin{itemize}
\item $\mathcal{S}$ has a right unit
$e_{\mathcal{S}}$, i.e.,
\begin{equation}\label{eq1.5.1} 
fe_{\mathcal{S}}=f\quad\mbox{for $f\in\mathcal{S}$}.
\end{equation}
\item Every element $f$ of $\mathcal{S}$ has a left
inverse $f^{\stackrel{\ell}{-}1}$ with respect to the
right unit $e_{\mathcal{S}}$, i.e.,
\begin{equation}\label{eq1.5.2} 
f^{\stackrel{\ell}{-}1}f=e_{\mathcal{S}}.
\end{equation}
\end{itemize}
\item $\phi : \mathcal{S}\to \mathcal{G}$ is a map
satisfying
\begin{equation}\label{eq1.5.3} 
\phi(fg)=\phi(f)\phi(g), \quad
\phi(\mathcal{S})\mathcal{S}\subseteq\mathcal{S},
\quad \phi(e_{\mathcal{S}})f=f, \quad 
e_{\mathcal{S}}f=\phi(f)e_{\mathcal{S}}
\end{equation}
\begin{equation}\label{eq1.5.4} 
\phi(f)f^{\stackrel{\ell}{-}1}=e_{\mathcal{S}},
\end{equation}
\begin{equation}\label{eq1.5.5} 
f\phi(g)=fg, \quad \phi(\phi(f)g)=\phi(f)\phi(g),
\end{equation}
where $f, g\in \mathcal{S}$.
\end{enumerate}

By Proposition~\ref{pr1.4.2},
Proposition~\ref{pr1.4.3} and
Proposition~\ref{pr1.4.4}, every digroup $G$ produces
a standard triple $(\Lr_G, \Ll_G, \phi)$ on $G$.

\bigskip
Now let $(\mathcal{S}, \mathcal{G}, \phi )$ be a
standard triple on $\Omega$ . Consider the set 
$$\mathcal{G}\times \mathcal{S}:=\{ \; (\alpha , f) \;
| \; \alpha\in\mathcal{G}, \; f\in\mathcal{S} \; \}.$$

Note that $\mathcal{G}\times \mathcal{S}$ can be
regarded as a subset of
$\mathcal{T}(\Omega\times\Omega)$ in the following
way:
$$ (\alpha , f)(x,y):=(\alpha (x) , f(y))
\quad\mbox{for $(x,y)\in \Omega\times\Omega$} .$$

\medskip
We define two binary operations $\xl$ and $\xr$ on
$\mathcal{G}\times \mathcal{S}$ by
\begin{eqnarray}
\label{eq1.5.6} (\alpha , f)\xl (\beta ,
g):&=&(\alpha\beta , fg),\\
\label{eq1.5.7} (\alpha , f)\xr (\beta ,
g):&=&(\alpha\beta , \phi (f)g),
\end{eqnarray}
where $\alpha$, $\beta \in \mathcal{G}$ and $f$, $g\in
\mathcal{S}$.

\medskip
One can check that the
$\xl$ and $\xr$ satisfy the diassociative law, $(1,
e_{\mathcal{S}})$ satisfies (\ref{eq0.1.4}) and
(\ref{eq0.1.5}), and $(\alpha ^{-1},
f^{\stackrel{\ell}{-}1})$ is the Liu inverse of
$(\alpha , f)$ with respect to the identity $(1,
e_{\mathcal{S}})$. This proves that $\mathcal{G}\times
\mathcal{S}$ is a digroup with the identity $(1,
e_{\mathcal{S}})$.

Therefore, we have the following proposition.

\begin{proposition}\label{pr1.4.6} Every standard
triple on a set $\Omega$ produces a digroup which is a
subset of $\mathcal{T}(\Omega\times\Omega)$.
$\qquad\|$ \end{proposition}

\section{Three remarks about digroups}

We finish this paper with three remarks about
digroups.

\medskip
1. It is challenging and interesting to develop the
counterpart of group theory in the context of
digroups. Even some very simple concepts in group
theory may resist our effort. One of the examples is
the order of an element of a group. Although different
criteria can be used to judge whether a generalization
of group theory is satisfactory, we believe that a
satisfying generalization of group theory should
contain a satisfying counterpart of the order of an
element of a group. Finding a suitable counterpart of
the order of an element of a group is time-consuming,
but the solution we have got seems to be a satisfying
solution (\cite{Liu1}).

\medskip
2. Digroups arise automatically from our study of
generalizing the Lie correspondence between linear Lie
groups and linear Lie algebras. However, the
non-associative algebra objects corresponding to
digroups are not Leibniz algebras, and the group-like
objects corresponding to Leibniz algebras are not
digroups. The following table gives a simple
sketch-map of the possible generalizations of the Lie
correspondence.

\bigskip
\begin{tabular}{|lcl|}\hline
\multicolumn{3}{|c|}{\rule[-3mm]{0mm}{9mm}\bfseries
The possible generalizations of the Lie
correspondence:}\\ 
\qquad\quad $G_0=\{\, Groups
\,\}$&$\longleftrightarrow$& $\mathcal{L}_0=\{\, Lie
\;\, algebras \,\}$\\
\hline\hline
The group-like objects&&The generalizations of Lie
algebras\\ \hline\hline
$G_1=\{\, ??? \,\}$&$\longleftrightarrow$&
$\mathcal{L}_1=\{\, Leibniz \;\, algebras \,\}$\\
\hline
$G_2=\{\, Digroups \,\}$&$\longleftrightarrow$&
$\mathcal{L}_2=\{\, ??? \,\}$\\ \hline
$\begin{array}{c} \star\\ \star\\ \star \end{array}
$&$\longleftrightarrow$& $\begin{array}{c} \bullet \\
\bullet \\ \bullet \end{array} $\\ \hline
\end{tabular}  

\bigskip
In the correspondence above, $G_1$ is a class of
group-like objects which is different from digroups,
and the notion of the Liu inverse is indispensable in
the description of $G_1$. $\mathcal{L}_2$ is a class
of generalization of Lie algebras which is different
from Leibniz algebras, and the notion of the Liu
inverse is also indispensable in the description of
$\mathcal{L}_2$. The star part $\begin{array}{c}
\star\\ \star\\ \star \end{array} $ is not empty and
consists of other kinds of group-like objects. The
corresponding bullet part $\begin{array}{c} \bullet \\
\bullet \\ \bullet \end{array} $ is also not empty and
consists of other kinds of generalizations of Lie
algebras. However, quantum groups do not belong to the
star part and super Lie algebras do not belong the
bullet part. The definitions of $G_1$ and
$\mathcal{L}_2$ will be given in \cite{Liu2}

\medskip
3. Another application of the notion of the Liu
inverse is in the search for the counterpart of
commutative rings. Using the notion of the Liu
inverse, we have introduced a class of ring-like
objects (\cite{Liu3}). A typical example of the class
of ring-like objects is the set of endomorphisms of a
commutative digroup. Some fundamental notions in
commutative rings like integral domains, prime ideals
and fields have natural counterparts in the class of
ring-like objects. Hence, there seems to be a good
foundation in the class of ring-like objects to
reconsider the theory of commutative rings.

\bigskip
\bigskip
{\bf Acknowledgment} I thank NSERC to provide the
partial support for me from March 1999 to March 2003.

\end{document}